\documentstyle[12pt]{article}

\font\tenmsb=msbm10
\font\sevenmsb=msbm7
\font\fivemsb=msbm5
\newfam\msbfam
\textfont\msbfam=\tenmsb
\scriptfont\msbfam=\sevenmsb
\scriptscriptfont\msbfam=\fivemsb
\def\Bbb#1{\fam\msbfam\relax#1}

\renewcommand{\P}{{\Bbb P}}
\newcommand{\ke}{{\cal E}}
\newcommand{\kf}{{\cal F}}
\newcommand{\kg}{{\cal G}}

\newcommand{\kk}{{\cal K}}

\newcommand{\kn}{{\cal N}}
\newcommand{\ko}{{\cal O}}

\newcommand{\kq}{{\cal Q}}

\newcommand{\kt}{{\cal T}}

\begin{document}
\newtheorem{lemma}{Lemma}[section]
\newtheorem{proposition}[lemma]{Proposition}
\newtheorem{remark}[lemma]{Remark}
\newtheorem{example}[lemma]{Example}
\newtheorem{theorem}[lemma]{Theorem}
\newtheorem{definition}[lemma]{Definition}
\newtheorem{corollary}[lemma]{Corollary}
\newtheorem{sub}[lemma]{}

\title{{\Large\bf The Tangent Space at a Special Symplectic Instanton Bundle
on ${\bf{\rm I\hspace{-0.15em}P}_{2n+1}}$}}
\author{{\bf
Giorgio Ottaviani}
\and
{\bf G\"{u}nther Trautmann
}}
\date{\small October 1993 \hspace{3.0cm} alg-geom/9402005}
\maketitle
\thispagestyle{empty}
\vspace{1.0cm}

\tableofcontents
\vspace{1cm}

\thispagestyle{empty}
\setcounter{page}{0}
\newpage

\section*{Introduction}\addcontentsline{toc}{section}{Introduction}

Mathematical instanton bundles on $\P_3$ have their analogues in rank--$2n$
instanton bundles on odd dimensional projective spaces $\P_{2n+1}$.   The
families of special instanton bundles on these spaces, which generalize the
special 'tHooft bundles on $\P_3$, were constructed and described in \cite{OS}
and \cite{ST}.   More general instanton bundles have recently been constructed
in \cite{AO2}.   Let $MI_{2n+1}(k)$ denote the moduli space of all instanton
bundles on $\P_{2n+1}$ with second Chern class $c_2=k$.   In order to obtain a
first impression of this space it is important to know its tangent dimension
$h^1End(\ke)$ at a stable bundle $\ke$ and the dimension $h^2End(\ke)$ of the
space of obstructions to smoothness.

In this paper we prove that for a special symplectic
bundle $\ke \in MI_{2n+1}(k)$
\[
h^2 End(\ke) = (k-2)^2 {2n-1\choose 2}.
\]
Such bundles are stable by \cite{AO1}.
So for $n \ge 2$ the situation is quite different to that of $\P_3$, where
this number becomes zero, which was shown in \cite{HN}.   Since $H^iEnd(\ke) =
0$ for $i \ge 3$, our result and the Hirzebruch--Riemann--Roch formula, see
Remark \ref{2.4},
\[
h^1End(\ke) - h^2End(\ke) = -k^2{2n-1\choose 2} + k(8n^2) + 1-4n^2
\]
give
\[
h^1End(\ke) = 4(3n-1) k + (2n-5)(2n-1).
\]

Therefore the dimension of $MI_{2n+1}(k)$ grows linearly in $k$, whereas the
difference $h^1End(\ke) - h^2End(\ke)$ becomes negative for $n \ge 2$ and
grows quadratically in $k$.   A more important consequence, however, is that
in general $MI_{2n+1}(k)$ cannot be smooth at special symplectic bundles, see
section 4 and \cite{AO2}.

In order to derive our result we fix a 2--dimensional vector space $U$ and
consider the natural action of $SL(2)$ on $\P_{2n+1} = \P(U \otimes S^nU)$ as
in \cite{ST}.   The special instanton bundles are related to the
$SL(2)$--homomorphisms $\beta$, see \ref{1.4}, and are $SL(2)$--invariant.   We
prove that there is an isomorphism of $SL(2)$--representations
\[
    H^2(End\,\ke) \cong S^{k-3}(U) \otimes S^{k-3}(U) \otimes S^2(U \otimes
    S^{n-2}U).
\]
\vfill

{\bf Acknowledgement}\quad  This work was supported by the
{\bf Deutsche Forschungs
Gemeinschaft}.   The first author wishes to thank the Fachbereich Mathematik of
the University of Kaiserslautern, where this work was begun, for its
hospitality.
\newpage

\section*{Notation}

\begin{itemize}
\item Throughout the paper $K$ denotes an algebraically closed ground field of
characteristic 0.

\item $U$ denotes a 2--dimensional $K$--vector space, $S_n = S^nU$ its $n$th
symmetric power and $V_n = U \otimes S_n$.

\item There is the natural exact squence of $GL(U)$--equivariant maps for any
$k,\, n \ge 1$
\[
    0 \to \Lambda^2 U \otimes S_{k-1} \otimes S_{n-1} \buildrel
    \beta\over\to S_k
    \otimes S_n \buildrel \mu\over\to S_{k+n} \to 0
\]
where $\mu$ is the multiplication map and $\beta$ is defined by $(s\wedge t)
\otimes f \otimes g \mapsto sf \otimes tg - tf \otimes sg$.   This sequence
splits and leads to the Clebsch--Gordan decomposition of $S_k \otimes S_n$ by
induction.   When we tensorize the sequence with $U$ we obtain the exact
sequence
\[
    0 \to \Lambda^2 U \otimes S_{k-1} \otimes V_{n-1}
    \buildrel\beta\over\to S_k
    \otimes V_n \buildrel\mu\over\to V_{k+n} \to 0.
\]

\item $\P = \P_{2n+1} = \P V_n$ is the projective space of one dimensional
subspaces of $V_n$.

\item The terms vector bundle and locally free sheaf are used synonymously.

\item $\ko(d)$ denotes the invertible sheaf of degree $d$ on $\P$, $\Omega^p$
the locally free sheaf of differential $p$--forms on $\P$, such that
$\Omega^p(p) = \Lambda^p \kq^\vee$ where $\kq = \kt(-1)$ is the canonical
quotient bundle on $\P$.

\item We use the abbreviations $\kf(d) = \kf \otimes_{\ko} \ko(d)$ for any
sheaf $\kf$ of $\ko$--modules on $\P,\; H^i\kf = H^i(\kf) = H^i(\P, \kf),\;
h^i\kf = \dim\, H^i\kf$.   If $E$ is a finite dimensional $K$--vector space,
$E \otimes \ko$ denotes the sheaf of sections of the trivial bundle $\P \times
E$, and $E \otimes \kf = (E \otimes \ko) \otimes_\ko \kf$.   We also write
$m\kf = K^m \otimes \kf$.

\item We use the Euler sequence $0\to \Omega^1(1) \to V^\vee \otimes \ko \to
\ko (1) \to 0$ and the derived sequences in its Koszul complex $0 \to
\Omega^p(p) \to \Lambda^p V^\vee \otimes \ko \to \Omega^{p-1} (p) \to 0$
without extra mentioning.

\item $Ext^i(\kf, \kg) = Ext^i_\ko (\P, \kf, \kg)$ for any two $\ko$--modules
$\kf$ and $\kg$.
\end{itemize}
\newpage

\section{Instanton bundles}

\begin{sub}\label{1.1}
{\rm An instanton bundle on $\P = \P_{2n+1}$ with
instanton number $k$ or a $k$--instanton is an algebraic vector bundle $\ke$
on $\P$ satisfying:

    \begin{itemize}
    \item[(i)] $\ke$ has rank $2n$ and Chern polynomial $c(\ke) =
    (1-h^2)^{-k} = 1 + k h^2 + \ldots$.
    \item[(ii)] $\ke$ has natural cohomology in the range $-2n-1 \le d \le
    0$, that is for any $d$ in that range $h^i\ke(d) \not= 0$ for at most one
    $i$.
    \end{itemize}

A $k$--instanton bundle $\ke$ is called {\bf symplectic} if there is an
isomorphism $\ke \buildrel \varphi\over\to \ke^\vee$ satisfying $\varphi^\vee =
-\varphi$.   In this case the spaces $A$ and $B$ below are Serre--duals of
each other, since
$H^{2n}(\ke(-2n-1))^\vee \cong H^1\ke^\vee(-1) \cong H^1\ke(-1)$.

{\bf Remark}:  In the original definition in \cite{OS} the additional
conditions

\begin{itemize}
\item[(iii)] $\ke$ is simple, that is $Hom(\ke,\ke) = K$,
\item[(iv)] the restriction of $\ke$ to a general line is trivial
\end{itemize}

are imposed.   It was shown in \cite{AO1} that (iii) is already a consequence
of (i) and (ii).   Condition (iv) seems to be independent but we do not need
it in this paper.   By \cite{ST} special instantons satisfy (iv).}
\end{sub}

\begin{sub}\label{1.2}
{\rm Let now $A, B, C$ be vector spaces of dimensions $k, k,
2n(k-1)$ respectively.   A pair of linear maps
\[
    A \buildrel a \over \longrightarrow B \otimes \Lambda^2 V^\vee_n,\quad B
    \otimes V^\vee_n \buildrel b\over\longrightarrow C
\]
corresponds to a pair of sheaf homomorphisms
\[
    A \otimes \ko(-1) \buildrel\tilde{a}\over\longrightarrow B \otimes \Omega^1
    (1), \quad B \otimes \Omega^1(1) \buildrel\tilde{b}\over\to C \otimes \ko.
\]

Here $\tilde{a}$ is the composition of the induced homomorphisms $A \otimes
\ko(-1) \to B \otimes \Lambda^2 V^\vee_n \otimes \ko(-1) \longrightarrow B
\otimes \Omega^1(1)$ and $\tilde{b}$ is the composition of the induced
homomorphismus $B \otimes \Omega^1(1) \longrightarrow B \otimes V^\vee_n
\otimes
\ko \to C \otimes \ko$.   Conversely, $a$ and $b$ are determined by
$\tilde{a}$ and $\tilde{b}$ respectively as $H^0 (\tilde{a}(1))$ and
$H^0(\tilde{b}^\vee)^\vee$.   Moreover, the sequence
\begin{equation}
A \otimes \ko(-1) \buildrel \tilde{a}\over\longrightarrow B \otimes
\Omega^1(1) \buildrel \tilde{b}\over\longrightarrow C \otimes \ko
\end{equation}
is a complex if and only if the induced sequence
\[
    A \longrightarrow B \otimes \Lambda^2 V^\vee_n \longrightarrow C \otimes
    V^\vee_n
\]
is a complex.   We say that (1) is a {\bf monad} if it is a complex and if in
addition $\tilde{a}$ is a subbundle and $\tilde{b}$ is surjective.}
\end{sub}

\begin{proposition}\label{1.3}
The cohomology sheaf $\ke = Ker\,
\tilde{b}/Im\, \tilde{a}$ of a monad (1) is a $k$--instanton and conversely
any $k$--instanton is the cohomology of a monad (1).   Moreover, the spaces
$A, B, C$ of such a monad can be identified with \mbox{$H^{2n}\ke(-2n-1)$,}
$H^1\ke(-1),\; H^1\ke$ respectively.
\end{proposition}

Sketch of a proof:  if a monad (1) is given it is easy to derive the
properties of the definition.   Conversely using Beilinson's spectral
sequence, Riemann--Roch and in particular (ii), one obtains a monad with the
identification of the vector spaces as in \cite{OS}.   The map $b$ is then
nothing but the
natural map $H^1\ke(-1) \otimes V^\vee_n \to H^1\ke$ and the map $a$ is given
as the composition of the cup product
\[
    H^{2n}\ke(-2n-1) \otimes \Lambda^2V_n \to H^{2n}\ke \otimes \Omega^{2n-1}
    (-1)
\]
and the natural isomorphisms
\[
    H^{2n}\ke\otimes\Omega^{2n-1}(-1) \cong H^{2n-1}\ke \otimes
\Omega^{2n-2}(-1)
    \cong \ldots \cong H^1\ke(-1)
\]
arising from the Koszul sequences and condition (ii), see \cite{V} in case of
$\P_3$.\\

\begin{sub}\label{1.4}{\rm{\bf Existence and special instanton bundles}:
Using the special
structure $V_n = U \otimes S^n U$ and the Clebsch--Gordan type exact sequence
\[
    0 \longrightarrow \Lambda^2 U \otimes S_{k-2} \otimes V_{n-1}
    \buildrel\beta\over\longrightarrow S_{k-1} \otimes V_n
    \buildrel\mu\over\longrightarrow V_{k+n+1} \longrightarrow 0,
\]
see notation, we define the special homomorphism
\[
    S^\vee_{k-1} \otimes \Omega^1(1) \buildrel\tilde{b}\over\longrightarrow
    \Lambda^2 U^\vee \otimes S^\vee_{k-2} \otimes V^\vee_{n-1} \otimes \ko
\]
by $b = \beta^\vee$.   We denote $\kk = Ker\;(\tilde{b})$.   It was shown in
\cite{ST} that $\tilde{b}$ is surjective and that
\[
    H^0 \kk(1) \subset S^\vee_{k-1} \otimes H^0 \Omega^1(2)
\]
can be identified with a canonical injective $GL(U)$--homomorphism
\[
    S^\vee_{2n+k-1} \otimes \Lambda^2 U^\vee \buildrel\kappa\over\to
S^\vee_{k-1} \otimes
    \Lambda^2 V^\vee_n,
\]
dual to the map
\[
    S_{k-1} \otimes \Lambda^2V_n \to S_{2n+k-1} \otimes \Lambda^2U
\]
which is defined by $f \otimes (s \otimes g) \wedge (t \otimes h) \mapsto
(fgh) \otimes
(s \wedge t)$.

In order to construct instanton bundles we have to find $k$--dimensional
subspaces
\[
    A \subset S^\vee_{2n+k-1} \otimes \Lambda^2 U^\vee \subset
    S^\vee_{k-1}
    \otimes \Lambda^2 V^\vee_n
\]
such that the induced homomorphism $\tilde{a}$ is a subbundle.   By \cite{ST},
Lemma 3.7.1, this is the case exactly when $\P A \subset \P (S^\vee_{2n+k-1})$
does not meet the secant variety $Sec_n(C_{2n+k-1})$ of $(n-1)$--dimensional
secant planes of the canonical rational curve $C_{2n+k-1}$ of $\P
S^\vee_{2n+k-1}$, given by $u \mapsto u^{2n+k-1}$.   By dimension reasons such
subspaces exist, \cite{ST}, 3.7, and hence instanton bundles exist.

A $k$--instanton bundle $\ke$ is called {\bf special} if the map $b$ of its
monad is isomorphic to the $GL(U)$--homomorphism $\beta^\vee$, that is if
there are isomorphisms $\varphi$ and $\psi$ and $g \in GL(V_n)$ with the
commutative diagram

\unitlength1cm
\begin{picture}(9,4)
\put(3.75,1.5){$S^\vee_{k-1} \otimes V^\vee_n$}
\put(6.25,1.6){\vector(1,0){1}}
\put(6.75,1.7){$\beta^\vee$}
\put(7.75,1.5){$\Lambda^2 U^\vee \otimes S^\vee_{k-2} \otimes V^\vee_{n-1}.$}
\put(3.45,2.5){$\varphi \otimes g^\vee$}
\put(4.75,3.0){\vector(0,-1){1}}
\put(4.95,2.5){$\wr\wr$}
\put(8.65,2.5){$\psi$}
\put(8.95,3.0){\vector(0,-1){1}}
\put(9.05,2.5){$\wr\wr$}
\put(3.05,3.3){$H^1\ke(-1) \otimes V^\vee_n$}
\put(6.25,3.4){\vector(1,0){1}}
\put(6.75,3.5){$b$}
\put(8.45,3.3){$H^1\ke$}
\end{picture}

Whereas in \cite{ST} the family of all special $k$--instanton bundles was
described, examples of different types of general instanton bundles were found
in \cite{AO2}.}
\end{sub}

\begin{remark}\label{1.5}
{\rm If $\ke$ is special and symplectic then, in addition to
the special $GL(U)$--homomorphism $b = \beta^\vee$ of its monad, the map $a$
is given by an element $\alpha \in S^\vee_{2n+2k-2}$ as $a = \kappa \circ
\tilde{\alpha}$ where $S_{k-1} \buildrel \tilde{\alpha}\over\to
S^\vee_{2n+k-1}$ is defined by $\tilde{\alpha}(f)(g) = \alpha(fg)$ and
$S^\vee_{2n+k-1} \buildrel \kappa\over\to S^\vee_{k-1} \otimes \Lambda^2
V^\vee_n$ is as above, \cite{ST}, 4.3 and 5.8.   In particular $a$ is a
$GL(U)$--homomorphism, too, and can be represented by a persymmetric matrix.}
\end{remark}

\begin{remark}\label{1.6}
{\rm It is shown in \cite{AO1} that special symplectic
instanton bundles are stable in the sense of Mumford--Takemoto.}
\end{remark}
\newpage

\section{Representing $\bf Ext^2(\ke,\ke)$}

\begin{proposition}\label{2.1}
Let $\ke$ be a symplectic $k$--instanton and let $\kn$ be the kernel of the
monad (1).   Then $Ext^2(\ke,\ke) \cong H^2(\kn \otimes \kn)$.
\end{proposition}

Proof:\quad The monad (1) gives rise to the exact sequences
\begin{equation}
    0 \longrightarrow \kn \longrightarrow B \otimes \Omega^1(1) \buildrel
    \tilde{b}\over\longrightarrow C \otimes \ko \longrightarrow 0
\end{equation}
and
\begin{equation}
    0 \longrightarrow A \otimes \ko (-1) \longrightarrow \kn \longrightarrow
\ke
    \longrightarrow 0.
\end{equation}

After tensoring we have the exact sequences
\begin{equation}
    0 \longrightarrow A \otimes \kn(-1) \longrightarrow \kn \otimes \kn
    \longrightarrow \ke \otimes \kn \longrightarrow 0
\end{equation}
and
\begin{equation}
    0 \longrightarrow A \otimes \ke (-1) \longrightarrow \kn \otimes \ke
    \longrightarrow \ke \otimes \ke \longrightarrow 0.
\end{equation}

Since $\ke \cong \ke^\vee$ we obtain $Ext^2(\ke,\ke) \cong H^2(\ke \otimes
\ke)$.   Sequence (2) implies $h^2\kn(-1) = h^3\kn(-1) = 0$ and from this and
(3) also $h^2\ke(-1) = h^3\ke(-1) = 0$.   Now sequences (4) and (5) yield
isomorphisms $H^2(\ke \otimes \ke) \cong H^2(\kn \otimes \ke) \cong H^2(\kn
\otimes \kn)$.\hfill$\Box$

\begin{sub}\label{2.2}
{\rm In order to represent $H^2(\kn \otimes \kn)$ we note that the sequence
(2) is part of the exact diagram
\begin{equation}
\begin{array}{ccccccccc}
  &      &     0      &     &    0                &    &    &    &\\
  &      & \downarrow &     & \downarrow &        &    &    &\\
0 & \longrightarrow  &    \kn     & \longrightarrow & B\otimes\Omega^1(1) &
\buildrel\tilde{b}\over\longrightarrow
& C \otimes \ko & \longrightarrow & 0\\
  &      & \downarrow &     & \downarrow &        &     \| &    &\\
0 & \longrightarrow  & H \otimes\ko& \longrightarrow &B\otimes
V^\vee\otimes \ko & \buildrel b\over\longrightarrow
& C \otimes \ko & \longrightarrow & 0\\
  &      & \downarrow & & \downarrow & & & & \\
  &      & B \otimes \ko(1) & = & B \otimes \ko(1) & & & & \\
  &      & \downarrow  & & \downarrow & & & &\\
  &      & 0 & & 0 & & & &
\end{array}
\end{equation}
where $H$ is the kernel of the operator $b$, which is surjective because
$\tilde{b}$ is surjective.   The
left--hand column of (6) gives us after tensoring by $\Omega^1(1)$
\begin{equation}
    B \otimes H^0\Omega^1(2) \buildrel\delta\over\cong H^1(\kn \otimes
    \Omega^1(1)) \mbox{ and } H^2(\kn \otimes \Omega^1(1)) = 0.
\end{equation}

Since $\tilde{b}$ is the Beilinson representation of $\kn$, we have the
commutative diagram
\begin{equation}
\begin{array}{ccc}
H^1\kn(-1) \otimes H^0 \ko(1) & \buildrel cup\over\longrightarrow &
H^1\kn\\
\|\wr & & \|\wr\\
B \otimes V^\vee & \buildrel b\over\longrightarrow & C.
\end{array}
\end{equation}

Moreover, $\delta$ in (7) coincides also with cup:
\begin{equation}
\begin{array}{cl}
B \otimes H^0 \Omega^1(2) &
\buildrel\delta\over{\buildrel\longrightarrow\over\approx}  H^1(\kn \otimes
\Omega^1(1))\\
\|\wr & \nearrow_{\mbox{cup}}\\
H^1\kn(-1)\otimes H^0\Omega^1(2)  &
\end{array}
\end{equation}

Tensoring the top row of (6) with $\kn$ and using (7) we obtain the following
diagram with exact row:
\begin{equation}
\begin{array}{ccccc}
0 \to H^1(\kn\otimes\kn)\to
& B \otimes H^1(\kn \otimes \Omega^1(1))
& \to
& C \otimes H^1(\kn)
& \to H^2(\kn \otimes\kn)  \to 0\\
& \|\wr & & \|\wr & \\
& B \otimes B \otimes \Lambda^2 V^\vee_n & \buildrel\Phi\over\longrightarrow
& C \otimes C. &
\end{array}
\end{equation}

It follows that
\begin{equation}
H^2(\kn \otimes \kn) = Coker(\Phi) = Ker(\Phi^\vee)^\vee.
\end{equation}
}\end{sub}

\begin{lemma}\label{2.3}
The induced operator $\Phi$ is the composition $B \otimes B \otimes \Lambda^2
V^\vee_n \buildrel id \otimes \sigma\over\longrightarrow B \otimes B \otimes
V^\vee_n \otimes V^\vee_n \buildrel b \otimes b\over\longrightarrow C \otimes
C$, where $\sigma$ denotes the canonical desymmetrization.
\end{lemma}

Proof:\quad  The computation of $\Phi$ is achieved by the diagram

\unitlength1cm
\begin{picture}(15,9)
\put(1,8){$B\otimes B\otimes \wedge^2 V^\vee_n$}
\put(0,6){$B\otimes H^1 \kn (-1) \otimes H^0\Omega^1(2)$}
\put(0.50,4){$B\otimes H^1(\kn\otimes \Omega^1(1))$}
\put(1.2,2){$C\otimes H^1\kn$}
\put(1.5,0){$C\otimes C$}

\put(2.1,7.6){\vector(0,-1){1}}
\put(2.1,5.6){\vector(0,-1){1}}
\put(2.1,3.6){\vector(0,-1){1}}
\put(2.1,1.6){\vector(0,-1){1}}

\put(2.2,7.1){$\approx$}
\put(2.2,5.1){$id_B \otimes  cup$}
\put(2.2,3.1){$H^1 (id_\kn \otimes \tilde{b})$}
\put(2.2,1.1){$\approx$}

\put(5.2,8.1){\vector(1,0){4}}
\put(5.2,6.1){\vector(1,0){4}}
\put(5.2,4.1){\vector(1,0){4}}

\put(6.2,8.3){$id_{B\otimes B}\otimes \sigma$}
\put(6.2,6.3){$id\otimes H^0(\iota (1))$}
\put(5.7,4.3){$id_B\otimes H^1(id_\kn\otimes \iota)$}

\put(10.7,8){$B\otimes B\otimes V^\vee_n \otimes V^\vee_n$}
\put(9.7,6){$B\otimes H^1 \kn(-1)\otimes V^\vee_n\otimes H^0 \ko(1)$}
\put(10.7,4){$B\otimes V^\vee_n \otimes H^1 \kn$}
\put(10.9,2){$B\otimes V^\vee_n \otimes C$}

\put(12.4,7.6){\vector(0,-1){1}}
\put(12.4,5.6){\vector(0,-1){1}}
\put(12.4,3.6){\vector(0,-1){1}}

\put(12.6,7.1){$\approx$}
\put(12.6,5.1){$id_{B\otimes V^\vee_n} \otimes  cup$}
\put(12.6,3.1){$\approx$}

\put(9.1,1.6){\vector(-4,-1){5}}
\put(9.1,3.6){\vector(-4,-1){5}}

\put(6.7,2.5){$b\otimes id_{H^1\kn}$}
\put(6.7,0.5){$b\otimes id_C$}

\put(7.3,7.2){I}
\put(7.3,5.2){II}
\put(7.3,3.5){III}
\put(7.3,1.8){IV}

\end{picture}

In this diagram $\iota$ denotes the canonical inclusion $\Omega^1(1)
\hookrightarrow V^\vee_n \otimes \ko$, and up to $\Lambda^2 V^\vee_n \cong
H^0\Omega^1(2)$ and $V^\vee_n \cong H^0 \ko(1)$ the map $\sigma$ can be
identified with $H^0(\iota(1))$.   Therefore, the square I is commutative.
Square II is a canonically induced diagram of cup--operations and commutative
using $B \cong H^1 \kn(-1)$.   The triangle III is induced by the commutative
triangle
\[
\begin{array}{cl}
    B\otimes \kn \otimes \Omega^1(1) & \buildrel id\otimes
    \iota\over\longrightarrow  B \otimes V^\vee \otimes \kn\\
    \quad\downarrow \tilde{b} \otimes id & \swarrow_{b \otimes id}  \\
    C \otimes \kn  &
\end{array}
\]
and hence commutative, and the commutativity of IV results just from the
identification $H^1 \kn \cong C$.   Now by definition the composition of the
left--hand column is $\Phi$ and the composition of the right--hand column is
$id_B \otimes id_{V^\vee_n} \otimes b$ since $b$ is defined by (8).

It follows that $\Phi = (b \otimes id_C) \circ (id_B \otimes id_{V^\vee_n}
\otimes b) \circ (id_{B \otimes B} \otimes \sigma) = (b \otimes b) \circ (id
\otimes \sigma)$.

\begin{remark}\label{2.4}
{\rm If $\ke$ is a $k$--instanton bundle it is easily checked
that $h^i\ke(d) = h^i\ke^\vee(d) = 0$ for $i \ge 2$ and $d \ge -1$.   Using
$\ke^\vee \otimes \kn$ again it follows that $Ext^i(\ke,\ke)
= H^i(\ke^\vee\otimes \ke) = H^i(\ke^\vee \otimes \kn) = 0$ for  $i \ge 3$.
This and the Riemann--Roch formula, which can also ad hoc be derived from the
monad representation, give
\[
h^1(\ke^\vee \otimes \ke) - h^2(\ke^\vee\otimes \ke) = -k^2{2n-1\choose 2} +
8kn^2 - 4n^2 + 1.
\]}
\end{remark}
\newpage

\section{Determination of $\bf Ext^2(\ke,\ke)$}

We are now able to determine $Ext^2(\ke,\ke)$ as a $GL(2)$--representation
space in case of a special instanton bundle.   In that case $b$ is the dual of
the operator $\beta : \Lambda^2 U \otimes S_{k - 2} \otimes V_{n-1} \to
S_{k - 1} \otimes V_n$, see notation or \ref{1.4}.   Then $\Phi^\vee$ is
the composition of $\beta \otimes \beta$ and the multiplication map $V_n
\otimes V_n \to \Lambda^2V_n$.   In order to simplify we choose a fixed basis
$s,t \in U$ and the isomorphism $\Lambda^2 U \cong k$ given by $s \wedge t$.
Then
\[
    S_{k - 2} \otimes S_{k - 2} \otimes V_{n-1} \otimes V_{n-1}
    \buildrel \Phi^\vee\over\to S_{k - 1} \otimes S_{k - 1} \otimes
    \Lambda^2 V_n
\]
is explicitly given by
\begin{eqnarray*}
    \Phi^\vee(g \otimes g' \otimes v \otimes v') & = & sg \otimes sg'
\otimes (tv
    \wedge tv') - sg \otimes tg' \otimes (tv \wedge sv')\\
    & - & tg \otimes sg' \otimes (sv \wedge tv') + tg \otimes tg' \otimes (sv
    \wedge sv').
\end{eqnarray*}

In order to determine the kernel of $\Phi^\vee$ we consider the
$GL(U)$--homomorphism
\[
    S_{k-3} \otimes S_{k-3} \otimes V_{n-2} \otimes V_{n-2} \buildrel
    \epsilon'\over\to S_{k-2} \otimes S_{k-2} \otimes V_{n-1} \otimes V_{n-1}
\]
defined similarly by
\begin{eqnarray*}
    \epsilon' (f \otimes f' \otimes u \otimes u') & = & sf \otimes sf' \otimes
    tu \otimes tu' - sf \otimes tf' \otimes su \otimes tu'\\
    & - & tf \otimes sf' \otimes tu \otimes su' + tf \otimes tf' \otimes su
    \otimes su'.
\end{eqnarray*}

Up to the order of factors the map $\epsilon'$ is the tensor product $\beta'
\otimes \beta'$ where $\beta' : S_{k-3} \otimes V_{n-2} \to S_{k-2} \otimes
V_{n-1}$ is defined as $\beta$.   Hence, $\epsilon'$ is injective.   Finally,
we define $\epsilon$ as the composition
\[
    S_{k-3} \otimes S_{k-3} \otimes S^2V_{n-2} \stackrel{id\otimes\iota}{
    \longrightarrow} S_{k-3} \otimes S_{k-3} \otimes V_{n-2} \otimes
    V_{n-2} \buildrel\epsilon'\over\to S_{k-2} \otimes S_{k-2}
\otimes V_{n-1} \otimes V_{n-1}
\]
where $\iota$ is the canonical desymmetrization.   Then also $\epsilon$ is
injective.

\begin{proposition}\label{3.1}
$(S_{k-3} \otimes S_{k-3} \otimes S^2V_{n-2}, \epsilon)$ is the kernel of
$\Phi^\vee$.
\end{proposition}

Proof:  A straightforward computation shows that $Im(\epsilon) \subset
Ker(\Phi^\vee)$.   In order to show equality we reduce $Ker(\Phi^\vee)$ modulo
$Im(\epsilon)$ using canonical bases of the vector spaces.   A more elegant
proof using Clebsch--Gordan decompositions seems much harder to achieve.   Let
us denote the bases as follows:

\begin{tabular}{lll}
    basis of  $S_{k-3}:$ & $e_\alpha = s^{k-3-\alpha} t^\alpha$ &$ 0 \le \alpha
    \le k-3$\\[0.5ex]
    basis of  $S_{k-2}:$ & $f_\alpha = s^{k-2-\alpha}t^\alpha $&$ 0 \le \alpha
    \le k-2$\\[0.5ex]
    basis of  $S_{k-1}:$ &$ g_\alpha = s^{k-1-\alpha}t^\alpha $&$ 0 \le
    \alpha \le k-1$\\[1.0ex]
    basis of  $V_{n-2}: $&$ u_\mu = s \otimes s^{n-2-\mu}t^\mu$ &$ 0 \le \mu
    \le n-2$\\
    & $\bar{u}_\mu = t \otimes s^{n-2-\mu} t^\mu$ & \\[0.5ex]
    basis of  $V_{n-1}: $&$ x_\mu = s \otimes s^{n-1-\mu}t^\mu$ &$ 0 \le \mu
    \le n-1$\\
     & $\bar{x}_\mu = t \otimes s^{n-1-\mu} t^\mu $& \\[0.5ex]
    basis of  $V_n: $&$ y_\mu = s \otimes  s^{n-\mu} t^\mu$ &$ 0 \le \mu
    \le n$\\
    & $\bar{y}_\mu = t \otimes s^{n-\mu}t^\mu.$ &
                          \end{tabular}

For the basis $f_\alpha \otimes f_\beta \otimes x_\mu \otimes x_\nu,\;
f_\alpha \otimes f_\beta \otimes x_\mu \otimes \bar{x}_\nu,\; f_\alpha \otimes
f_\beta \otimes \bar{x}_\mu \otimes x_\nu,\, f_\alpha \otimes f_\beta \otimes
\bar{x}_\mu \otimes \bar{x}_\nu$ we use the index tuplets $(\alpha, \beta,
\mu, \nu),\; (\alpha, \beta, \mu,\bar{\nu}),\; (\alpha, \beta, \bar{\mu},\nu),
\; (\alpha, \beta, \bar{\mu}, \bar{\nu})$ respectively.   The set
of these indices will be ordered {\bf lexicographically} with the additional
assumption that always $\mu < \bar{\nu}$.   Then, for example, $(\alpha,
\beta, \mu, \bar{\nu}) < (\alpha, \beta, \bar{\lambda}, \delta)$.

Accordingly, the coefficients of an element $\xi \in S_{k-2} \otimes S_{k-2}
\otimes V_{n-1} \otimes V_{n-1}$ will be denoted by $c(\alpha, \beta, \mu,
\nu),\; c(\alpha, \beta, \mu, \bar{\nu}),\; c(\alpha, \beta, \bar{\mu},
\nu),\; c(\alpha, \beta, \bar{\mu}, \bar{\nu})$.

By the formula for $\Phi^\vee$ we obtain the

\begin{lemma}\label{3.2}
Let $\xi \in S_{k-2} \otimes S_{k-2} \otimes V_{n-1} \otimes V_{n-1}$.

\begin{itemize}
\item[(i)] The coefficient of $\Phi^\vee(\xi)$ at the basis element $g_\alpha
\otimes g_\beta \otimes y_\mu \wedge \bar{y}_\nu$ in $S_{k-1} \otimes S_{k-1}
\otimes \Lambda^2V_n$ is

$\phantom{-}c(\alpha, \beta, \mu-1, \overline{\nu-1}) - c(\alpha, \beta,
\overline{\nu-1}, \mu-1)$\\[0.5ex]
$-c(\alpha, \beta-1, \mu-1, \bar{\nu}) + c(\alpha, \beta-1, \overline{\nu-1},
\mu)$\\[0.5ex]
$-c(\alpha-1, \beta, \mu, \overline{\nu-1}) + c(\alpha-1, \beta, \bar{\nu},
\mu)$\\[0.5ex]
$+c(\alpha-1, \beta-1, \mu, \bar{\nu}) - c(\alpha-1, \beta-1, \bar{\nu},
\bar{\mu}).$

Here we agree that each of these coefficients is 0 if one of $\alpha,
\alpha-1, \beta,
\beta-1 \not\in [0,k-2]$ or if one of $\mu, \mu-1, \nu, \nu-1 \not\in [0,
n-1]$.

\item[(ii)] Analogous statements hold for the coefficient of $\Phi^\vee(\xi)$
at $g_\alpha \otimes g_\beta \otimes y_\mu \wedge y_\nu$ for $\mu < \nu$
(without bars) and at $g_\alpha \otimes g_\beta \otimes \bar{y}_\mu \wedge
\bar{y}_\nu$ for $\mu < \nu$ (with two bars).
\end{itemize}
\end{lemma}

\begin{lemma}\label{3.3}  Let the notation be as above.   If $\Phi^\vee(\xi) =
0$ then:

\begin{itemize}
\item[(i)] If $c(\alpha, \beta, \mu, \nu)$ is the first non--zero coefficient
of $\xi$ (in the lexicographical order), then $0 < \mu \le \nu$.

\item[(ii)] If $c(\alpha, \beta, \mu, \bar{\nu})$ is the first non--zero
coefficient of $\xi$, then $\mu \not= 0,\; \nu \not= 0$.

\item[(iii)] $c(\alpha, \beta, \bar{\mu}, \nu)$ is never a first non--zero
coefficient of $\xi$.

\item[(iv)] If $c(\alpha, \beta, \bar{\mu}, \bar{\nu})$ is the first non--zero
coefficient of $\xi$, then $0 < \mu \le \nu$.
\end{itemize}
\end{lemma}

Proof:  (i)\quad Let $c(\alpha, \beta, \mu, \nu)$ be the first coefficient of
$\xi$.   Then, by Lemma \ref{3.2} the coefficient of $0 = \Phi^\vee(\xi)$ at
$g_\alpha \otimes g_\beta \otimes y_{\mu+1} \wedge y_{\nu+1}$ is
\begin{eqnarray*}
    0 & = & c(\alpha, \beta, \mu, \nu) - c(\alpha, \beta, \nu, \mu) - c(\alpha,
    \beta-1, \mu, \nu + 1) + c(\alpha, \beta-1, \nu, \mu + 1)\\
      & - & c(\alpha-1, \beta, \mu + 1, \nu) +
c(\alpha-1, \beta, \nu + 1, \mu) -
      \ldots
\end{eqnarray*}

Since $c(\alpha, \beta, \mu, \nu)$ is the first coefficient, only the first two
in this formula could be non--zero because the others have smaller index in
the lexicographical order.   Hence
\[
    c(\alpha, \beta, \mu, \nu) = c(\alpha, \beta, \nu, \mu).
\]

If $\mu > \nu$ then $c(\alpha, \beta, \nu, \mu)$ would be earlier and
non--zero.   Hence, $\mu \le \nu$.   Assume now that $\mu = 0$.   The
coefficient of $\phi^\vee(\xi)$ of $g_\alpha \otimes g_{\beta + 1} \otimes y_0
\wedge y_{\nu+1}$ is
\begin{eqnarray*}
    0 & = & c(\alpha, \beta + 1, -1, \nu) - c(\alpha, \beta+1, \nu, -1)\\
      & - & c(\alpha, \beta, -1, \nu+1) +   c(\alpha, \beta, \nu, 0) \mp \ldots
\end{eqnarray*}

In this sum all but $c(\alpha, \beta, \nu, 0)$ are automatically zero because
$(\alpha-1, \beta, \ldots) \le (\alpha, \beta, 0, \nu)$ and $-1$ occurs.
Hence, $c(\alpha, \beta, 0, \nu) = c(\alpha, \beta, \nu, 0) = 0$,
contradiction.

The statements (ii), (iii), (iv) are proved analogously. \hfill $\Box$\\

Now we continue the proof of Proposition \ref{3.1}.   We reduce an element
$\xi \in Ker(\Phi^\vee)$ to $0\, mod\, Im(\epsilon)$ using  Lemma \ref{3.3}.

\begin{itemize}
\item[a)] Assume that the first non--zero coefficient of $\xi$ is
\[
    c(\alpha, \beta, \mu, \nu).
\]

Then by Lemma \ref{3.3}\quad  $0 < \mu \le \nu$.
Then the element
\[
    \xi' = \xi - c(\alpha, \beta, \mu, \nu) \epsilon (e_\alpha \otimes e_\beta
    \otimes u_{\mu-1} \cdot u_{\nu-1})
\]
belongs to $Ker(\Phi^\vee)$.   We have

$\epsilon(e_\alpha \otimes e_\beta \otimes u_{\mu-1} \cdot u_{\nu-1})$\\
$ = f_\alpha \otimes f_\beta \otimes (x_\mu \otimes x_\nu + x_\nu \otimes
x_\mu)$\\
$- f_\alpha \otimes f_{\beta+1} \otimes (x_{\mu-1} \otimes x_\nu + x_{\nu-1}
\otimes x_\mu)$\\
$- f_{\alpha + 1} \otimes f_\beta \otimes (x_\mu \otimes x_{\nu-1} + x_\nu
\otimes x_{\mu-1})$\\
$+ f_{\alpha +1} \otimes f_{\beta + 1} \otimes (x_{\mu-1} \otimes x_{\nu-1} +
x_{\nu-1} \otimes x_{\mu-1})$

and therefore $\xi'$ is a sum of monomials of index $> (\alpha, \beta, \mu,
\nu)$.   Hence, we can assume that $\xi\, mod\, Im(\epsilon)$ has no
coefficient with index $(\alpha, \beta, \mu, \nu)$.

\item[b)] By Lemma \ref{3.3} we can assume that the first non--zero
coefficient of $\xi$ has index $(\alpha, \beta, \mu, \bar{\nu})$ or $(\alpha,
\beta, \bar{\mu}, \bar{\nu})$.   In the first case we know by Lemma \ref{3.3}
that $0 < \mu, \nu$.   When we consider again
\[
    \xi' = \xi - c(\alpha, \beta, \mu, \bar{\nu}) \epsilon (e_\alpha \otimes
    e_\beta \otimes u_{\mu-1} \cdot \bar{u}_{\nu-1})
\]
we have $\phi^\vee(\xi') = 0$ and $\xi'$ is a sum of monomials of index
$>(\alpha, \beta, \mu, \bar{\nu})$.   Hence, we may assume that $\xi\, mod\,
Im(\epsilon)$ has $c(\alpha, \beta, \bar{\mu}, \bar{\nu})$ as first non--zero
coefficient.   Again by Lemma \ref{3.3}\quad $0 < \mu, \nu$ and
\[
    \xi' = \xi - c(\alpha, \beta, \bar{\mu}, \bar{\nu}) \epsilon (e_\alpha
    \otimes e_\beta \otimes \bar{u}_{\mu-1} \cdot \bar{u}_{\nu-1})
\]
is a sum of monomials of index $>(\alpha, \beta, \bar{\mu}, \bar{\nu})$.

This finally shows that $\xi = 0\, mod\, Im(\epsilon)$.
\end{itemize}

This completes the proof of Proposition \ref{3.1}.
\newpage

\section{Conclusions}

By Proposition \ref{2.1}, Proposition \ref{3.1}, (11) and Lemma \ref{2.3}
we have determined the space
$Ext^2(\ke, \ke)$.   Together with Remark \ref{2.4} we obtain

\begin{theorem}\label{4.1}  For any special symplectic $k$--instanton bundle
$\ke$ on $\P_{2n+1}$

\begin{itemize}
\item[(1)] $Ext^2(\ke,\ke) \cong S^\vee_{k-3} \otimes S^\vee_{k-3} \otimes
S^2V^\vee_{n-2}$

\item[(2)] $\dim\, Ext^2(\ke,\ke) = (k-2)^2 {2n-1\choose 2}$

\item[(3)] $\dim\, Ext^1(\ke,\ke) = 4k(3n-1) + (2n-5)(2n-1)$.
\end{itemize}
\end{theorem}

Let $MI_{2n+1}(k)$ denote the open part of the Maruyama scheme of
semi--stable coherent sheaves on $\P_{2n+1}$ with Chern polynomial
$(1-h^2)^{-k}$
consisting of instanton bundles.   By \cite{AO1} any special symplectic
instanton bundle $\ke$ is stable.   Therefore, $Ext^1(\ke,\ke)$ can be
identified with the tangent space of $MI_{2n+1}(k)$ at $\ke$.   In \cite{AO2}
deformations $\ke'$ of special symplectic instanton bundles in $MI_{2n+1}(k)$
have been found for $n = 2$ and $k = 3,4$ which satisfy $Ext^2(\ke', \ke') =
0$.
This shows that in these cases there
are components $MI'_{2n+1}(k)$ of $MI_{2n+1}(k)$ of the expected dimension
$4(3n-1) k+(2n-5)(2n-1)$ containing the set of special instanton bundles.   In
particular, see \cite{AO2}:

{\it for $k = 3,4$  the moduli space
$MI_5(k)$ is singular at least in special symplectic bundles.}

However, in case $2n+1 = 3$ we obtain the vanishing result of \cite{HN}:

{\it
any special $k$--instanton bundle $\ke$ on $\P_3$ satisfies $Ext^2(\ke,\ke) =
0$ and is a smooth point of $MI_3(k)$,}

since any rank--2 instanton bundle is
symplectic.

\vfill

\begin{tabular}{ll}
Giorgio Ottaviani                    & G\"{u}nther Trautmann\\
Dipartimento di Matematica Applicata & Fachbereich Mathematik\\
Via S.\ Marta 3                      & Erwin--Schr\"{o}dinger Stra{\ss}e\\
50139 Firenze                        & 67663 Kaiserslautern\\
Italy                                & Germany\\[1.0ex]
ottaviani@ingf\/i1.ing.unif\/i.it            & trm@mathematik.uni-kl.de
\end{tabular}
\end{document}